\documentclass[12pt,a4paper]{article}
\usepackage{amsmath}
\usepackage{amssymb}
\usepackage{t1enc}
\usepackage[latin1]{inputenc}
\usepackage[english]{babel}
\usepackage{enumitem}
\usepackage{amsfonts}
\usepackage{latexsym}
\usepackage{bm}
\newtheorem{theorem}{Theorem}[section]
\newtheorem{prop}[theorem]{Proposition}
\newtheorem{lemma}[theorem]{Lemma}

\newtheorem{conj}[theorem]{Conjecture}

\newtheorem{problem}[theorem]{Problem}

\title{Cross-intersecting non-empty uniform subfamilies of hereditary families}
\author{Peter Borg\\[5mm]
Department of Mathematics, University of Malta, Malta\\
\texttt{peter.borg@um.edu.mt}}
\date{}

\begin{document}
\maketitle

\begin{abstract}
A set $A$ \emph{$t$-intersects} a set $B$ if $A$ and $B$ have at least $t$ common elements. A set of sets is called a \emph{family}. Two families $\mathcal{A}$ and $\mathcal{B}$ are \emph{cross-$t$-intersecting} if each set in $\mathcal{A}$ $t$-intersects each set in $\mathcal{B}$.  A family $\mathcal{H}$ is \emph{hereditary} if for each set $A$ in $\mathcal{H}$, all the subsets of $A$ are in $\mathcal{H}$. The \emph{$r$th level of $\mathcal{H}$}, denoted by $\mathcal{H}^{(r)}$, is the family of $r$-element sets in $\mathcal{H}$. A set $B$ in $\mathcal{H}$ is a \emph{base of $\mathcal{H}$} if for each set $A$ in $\mathcal{H}$, $B$ is not a proper subset of $A$. Let $\mu(\mathcal{H})$ denote the size of a smallest base of $\mathcal{H}$. We show that for any integers $t$, $r$, and $s$ with $1 \leq t \leq r \leq s$, there exists an integer $c(r,s,t)$ such that the following holds for any hereditary family $\mathcal{H}$ with $\mu(\mathcal{H}) \geq c(r,s,t)$. If $\mathcal{A}$ is a non-empty subfamily of $\mathcal{H}^{(r)}$, $\mathcal{B}$ is a non-empty subfamily of $\mathcal{H}^{(s)}$, $\mathcal{A}$ and $\mathcal{B}$ are cross-$t$-intersecting, and $|\mathcal{A}| + |\mathcal{B}|$ is maximum under the given conditions, then for some set $I$ in $\mathcal{H}$ with $t \leq |I| \leq r$, either $\mathcal{A} = \{A \in \mathcal{H}^{(r)} \colon I \subseteq A\}$ and $\mathcal{B} = \{B \in \mathcal{H}^{(s)} \colon |B \cap I| \geq t\}$, or $r = s$, $t < |I|$, $\mathcal{A} = \{A \in \mathcal{H}^{(r)} \colon |A \cap I| \geq t\}$, and $\mathcal{B} = \{B \in \mathcal{H}^{(s)} \colon I \subseteq B\}$. This was conjectured by the author for $t=1$ and generalizes well-known results for the case where $\mathcal{H}$ is a power set. 
\end{abstract}

\section{Introduction}

\subsection{Basic definitions and notation}
%Before introducing the problems that we will deal with, we provide the main definitions and notation, which will be used throughout the paper.
%In this section, we first provide definitions and notation that will be %used throughout the paper, and then we introduce the problems that we %will deal with.

Unless otherwise stated, we shall use small letters such as $x$ to denote non-negative integers or elements of a set, capital letters such as $X$ to denote sets, and calligraphic letters such as $\mathcal{F}$ to denote \emph{families} (that is, sets whose members are sets themselves). Arbitrary sets and families are taken to be finite and may be the \emph{empty set} $\emptyset$. An \emph{$r$-element set} is a set of size $r$, that is, a set having exactly $r$ elements (also called members). The set of positive integers is denoted by
$\mathbb{N}$. For $m, n \in \mathbb{N}$, the set $\{i \in \mathbb{N} \colon m \leq i \leq n\}$ is denoted by $[m,n]$. We abbreviate $[1,n]$ to $[n]$, and we take $[0]$ to be $\emptyset$. For a set $X$, the \emph{power set of $X$} (that is, $\{A \colon A \subseteq X\}$) is denoted by $2^X$, and the family $\{A \subseteq X \colon |A| = r\}$ is denoted by $X \choose r$.

We say that a set $A$ \emph{$t$-intersects} a set $B$ if $A$ and $B$ have at least $t$ common elements. A family $\mathcal{A}$ is said to be \emph{$t$-intersecting} if for every $A, B \in \mathcal{A}$, $A$ $t$-intersects $B$. A $1$-intersecting family is also simply called an \emph{intersecting family}. A $t$-intersecting family $\mathcal{A}$ is said to be \emph{trivial} if its sets have at least $t$ common elements. For a family $\mathcal{F}$ and a $t$-element set $T$, the family $\{A \in \mathcal{F} \colon T \subseteq A\}$ is denoted by $\mathcal{F}(T)$ and called a \emph{$t$-star of $\mathcal{F}$}. Note that non-empty $t$-stars are trivial $t$-intersecting families. We say that $\mathcal{F}$ has the \emph{$t$-star property} if at least one of the largest $t$-intersecting subfamilies of $\mathcal{F}$ is a $t$-star of $\mathcal{F}$.

\subsection{Intersecting families}

One of the most popular endeavours in extremal set theory is that of determining the size or the structure of a largest $t$-intersecting subfamily of a given family $\mathcal{F}$. This originated in \cite{EKR}, which features the classical result referred to as the Erd\H os-Ko-Rado (EKR) Theorem. The EKR Theorem says that for $1 \leq t \leq r$ there exists an integer $n_0(r,t)$ such that for $n \geq n_0(r,t)$, the size of a largest $t$-intersecting subfamily of ${[n] \choose r}$ is ${n-t \choose r-t}$, meaning that ${[n] \choose r}$ has the $t$-star property. It also says that the smallest possible $n_0(r,1)$ is $2r$; among the various proofs of this fact (see \cite{EKR,Kat,HM,K,D,FF2}) there is a short one by Katona \cite{K}, introducing the elegant cycle method, and another one by Daykin \cite{D}, using the Kruskal-Katona Theorem \cite{Kr,Ka}. Note that ${[n] \choose r}$ itself is intersecting if $n < 2r$. The EKR Theorem inspired a sequence of results \cite{F_t1,W,FF,AK1} that culminated in the complete solution of the problem for $t$-intersecting subfamilies of ${[n] \choose r}$. The solution had been conjectured by Frankl \cite{F_t1}. It particularly tells us that the smallest possible $n_0(r,t)$ is $(t+1)(r-t+1)$; this was established by Frankl \cite{F_t1} %for $t \geq 15$, and subsequently by 
and Wilson \cite{W}. % for any $t$. 
Ahlswede and Khachatrian \cite{AK1} settled the case $n < (t+1)(r-t+1)$. The $t$-intersection problem for $2^{[n]}$ was solved by Katona \cite{Kat}. These are among the most prominent results in extremal set theory. The EKR Theorem inspired a wealth of results that establish how large a system of sets can be under certain intersection conditions; see \cite{DF,F,F2,HST,HT,Borg7,FTsurvey}.

A set $B$ in a family $\mathcal{F}$ is called a \emph{base of $\mathcal{F}$} if for each $A \in \mathcal{F}$, $B$ is not a proper subset of $A$. The size of a smallest base of $\mathcal{F}$ is denoted by $\mu(\mathcal{F})$. The family of $r$-element sets in $\mathcal{F}$ is denoted by $\mathcal{F}^{(r)}$ and called the \emph{$r$th level of $\mathcal{F}$}. 

A family $\mathcal{F}$ is said to be \emph{hereditary} if for each $A \in \mathcal{F}$, all the subsets of $A$ are members of $\mathcal{F}$. In the literature, a hereditary family is also called an \emph{ideal}, a \emph{downset}, and an \emph{abstract simplicial complex}. Hereditary families are important combinatorial objects that have attracted much attention. The various interesting examples include the family of \emph{independent sets} of a \emph{graph} or a \emph{matroid}. The power set is the simplest example. In fact, by definition, a family is hereditary if and only if it is a union of power sets. Note that if $X_1, \dots, X_k$ are the bases of a hereditary family $\mathcal{H}$, then $\mathcal{H} = 2^{X_1} \cup \dots \cup 2^{X_k}$.  

The most basic result on intersecting families, also proved in the seminal EKR paper \cite{EKR}, is that the hereditary family $2^{[n]}$ has the $1$-star property. One of the central conjectures in extremal set theory, due to Chv\'atal \cite{Chv}, is that every hereditary family $\mathcal{H}$ has the $1$-star property. 
%
%\begin{conj}[\cite{Chv}] Every hereditary family has the $1$-star property.
%\end{conj}
%
Several cases have been verified \cite{Chva, Sterboul, Schonheim, Miklos2, Miklos, Wang, Sn} (see also \cite{Chvatalsite}), many of which are captured by Snevily's result \cite{Sn} % The closest result so far is due to Snevily \cite{Sn} %who used a result of Berge \cite{...} 
(\cite{Borg4b} provides a generalization obtained by means of a self-contained alternative argument). For $t \geq 2$, the $t$-star property fails already for $\mathcal{H} = 2^{[n]}$ with $n \geq t+2$; the largest $t$-intersecting subfamilies of $2^{[n]}$ were determined by Katona \cite{Kat}. However, for levels of hereditary families, we have the following generalization of the Holroyd--Talbot Conjecture \cite[Conjecture~7]{HT}.
\begin{conj}[\cite{Borg}] \label{AK gen} If $1 \leq t \leq r$ and $\mathcal{H}$ is a hereditary family with $\mu(\mathcal{H}) \geq (t+1)(r-t+1)$, then $\mathcal{H}^{(r)}$ has the $t$-star property.
\end{conj}
Note that if $\mathcal{H} = 2^{[n]}$, then $\mathcal{H}^{(r)} = {[n] \choose r}$ and $\mu(\mathcal{H}) = n$. It follows by the above-mentioned results for ${[n] \choose r}$ that the conjecture is true for $\mathcal{H}=2^{[n]}$ and that the condition $\mu(\mathcal{H}) \geq (t+1)(r-t+1)$ cannot be improved. The author verified the conjecture for $\mu(\mathcal{H})$ sufficiently large depending only on $r$ and $t$.
\begin{theorem}[\cite{Borg}]\label{t int her} Conjecture~\ref{AK gen} is true if $\mu(\mathcal{H}) \geq (r-t){3r-2t-1 \choose t+1} + r$.
\end{theorem}
By \cite[Theorem~1.2 and Section~4.1]{Borgmaxprod}, Conjecture~\ref{AK gen} is also true if $\mu(\mathcal{H}) \geq (r-t)r{r \choose t} + r$.

\subsection{Cross-intersecting families}

A popular variant of the intersection problem described above is the cross-intersection problem.

Two families $\mathcal{A}$ and $\mathcal{B}$ are said to be \emph{cross-$t$-intersecting} if each set in $\mathcal{A}$ $t$-intersects each set in $\mathcal{B}$. Cross-$1$-intersecting families are also simply called \emph{cross-intersecting families}. 

For $t$-intersecting subfamilies of a given family $\mathcal{F}$,
the natural question to ask is how large they can be. For
cross-$t$-intersecting families, two natural parameters arise: the sum and the product of sizes of the cross-$t$-intersecting families. The problem of maximizing the sum or the product of sizes of cross-$t$-intersecting subfamilies of a given family $\mathcal{F}$ has been attracting much attention (many of the results to date are referenced in \cite{Borg8, Borgmaxprod, BorgJLMS}). 

In this paper, we are concerned with the sum problem for the case where, as in Theorem~\ref{t int her}, $\mathcal{F}$ is a level of a hereditary family, but we also address the problem where the cross-$t$-intersecting families come from different levels and are non-empty. Thus, it is convenient to introduce the following notation. For two families $\mathcal{F}$ and $\mathcal{G}$, let 
\begin{center} $C(\mathcal{F},\mathcal{G},t) = \{(\mathcal{A},\mathcal{B}) \colon \emptyset \neq \mathcal{A} \subseteq \mathcal{F}, \emptyset \neq \mathcal{B} \subseteq \mathcal{G}, \mbox{$\mathcal{A}$ and $\mathcal{B}$ are cross-$t$-intersecting}\},$
\end{center}
\begin{center}  
$m(\mathcal{F},\mathcal{G},t) = \max\{|\mathcal{A}| + |\mathcal{B}| \colon (\mathcal{A},\mathcal{B}) \in C(\mathcal{F},\mathcal{G},t)\},$
\end{center}
\begin{center}
$M(\mathcal{F},\mathcal{G},t) = \{(\mathcal{A},\mathcal{B}) \in C(\mathcal{F},\mathcal{G},t) \colon  |\mathcal{A}| + |\mathcal{B}| = m(\mathcal{F},\mathcal{G},t)\}$.
\end{center}
%
%$C(\mathcal{F},\mathcal{G},t)$ denote the set $\{(\mathcal{A},\mathcal{B}) \colon \emptyset \neq \mathcal{A} \subseteq \mathcal{F}, \emptyset \neq \mathcal{B} \subseteq \mathcal{G},$ $\mathcal{A}$ and $\mathcal{B}$ are cross-$t$-intersecting, and if $\mathcal{F} = \mathcal{G}$, then $|\mathcal{A}| \leq |\mathcal{B}|\}$, let $m(\mathcal{F},\mathcal{G},t)$ denote $\max\{|\mathcal{A}| + |\mathcal{B}| \colon (\mathcal{A},\mathcal{B}) \in C(\mathcal{F},\mathcal{G},t)\}$, and let $M(\mathcal{F},\mathcal{G},t)$ denote the set of pairs $(\mathcal{A}, \mathcal{B})$ in $C(\mathcal{F},\mathcal{G},t)$ such that $|\mathcal{A}| + |\mathcal{B}|$ is maximum, that is, $M(\mathcal{F},\mathcal{G},t) = \{(\mathcal{A},\mathcal{B}) \in C(\mathcal{F},\mathcal{G},t) \colon  |\mathcal{A}| + |\mathcal{B}| = m(\mathcal{F},\mathcal{G},t)\}$. 

Hilton and Milner \cite{HM} showed that if $\mathcal{A}$ and $\mathcal{B}$ are non-empty cross-intersecting subfamilies of ${[n] \choose r}$ with $1 \leq r \leq n/2$, then $|\mathcal{A}| + |\mathcal{B}| \leq {n \choose r} - {n - r \choose r} + 1$. Equality holds if $\mathcal{A}$ consists of $[r]$ only and $\mathcal{B}$ consists of all the sets in ${[n] \choose r}$ that intersect $[r]$. In other words, if $1 = t \leq r \leq n/2$ and $\mathcal{F} = \mathcal{G} = {[n] \choose r}$, then %$m(\mathcal{F}, \mathcal{G}, 1) = {n \choose r} - {n - r \choose r} + 1$ and 
$(\{[r]\},\{B \in \mathcal{G} \colon B \cap [r] \neq \emptyset\}) \in M(\mathcal{F}, \mathcal{G}, t)$. %Simpson \cite{Simpson} gave a streamlined proof. 
%Wang and Zhang \cite{WZ} provided a significant generalization; in particular, they determined $m(\mathcal{F}, \mathcal{G}, t)$ for $\mathcal{F} = \mathcal{G} = {[n] \choose r}$. %with $1 \leq t \leq r \leq n$. 
Frankl and Tokushige \cite{FT1} showed that the same holds in the more general case where $1 = t \leq r \leq s$, $n \geq r+s$, $\mathcal{F} = {[n] \choose r}$, and $\mathcal{G} = {[n] \choose s}$. Wang and Zhang \cite{WZ2} generalized this for $t \geq 1$. They proved that if $t < \min\{r,s\}$, $n \geq r+s-t+1$, ${n \choose r} \leq {n \choose s}$, $\mathcal{F} = {[n] \choose r}$, and $\mathcal{G} = {[n] \choose s}$, then $(\{[r]\}, \{B \in \mathcal{G} \colon |B \cap [r]| \geq t\}) \in M(\mathcal{F}, \mathcal{G}, t)$ 
%
%\begin{theorem}[\cite{WZ}] \label{WZthm} If $1 \leq t < \min\{r,s\}$, $n > r+s-t$, ${n \choose r} \leq {n \choose s}$, $\mathcal{F} = {[n] \choose r}$, and $\mathcal{G} = {[n] \choose s}$, then $(\{[r]\}, \{B \in \mathcal{G} \colon |B \cap [r]| \geq t\}) \in M(\mathcal{F}, \mathcal{G}, t)$.
%\end{theorem}
%
(an independent proof for $r=s$ has been obtained by Frankl and Kupavskii \cite{FK}); they also determined the pairs in $M(\mathcal{F}, \mathcal{G}, t)$. %extremal structures. 
%This implies that $m(\mathcal{F}, \mathcal{G},t) = ...$ and %from Theorem~\ref{WZthm} 
It immediately follows that if we allow the cross-$t$-intersecting families $\mathcal{A}$ and $\mathcal{B}$ to be empty, then $|\mathcal{A}| + |\mathcal{B}|$ is maximum if $\mathcal{A} = \emptyset$ and $\mathcal{B} = {[n] \choose s}$. %; see Lemma~\ref{nonemptylemma} with $\mathcal{H} = 2^{{n}}$

%Wang and Zhang \cite{WZ} completely solved the maximum sum cross-$t$-intersection problem for ${[n] \choose r}$ by reducing it to the complete $t$-intersection theorem in \cite{AK1}, using an elegant combination of the method used in \cite{Borg4,Borg3,Borg2,Borg5} and an important lemma that is found in \cite{AC,CK} and referred to as the `no-homomorphism' lemma. The solution for the case $t = 1$ had been obtained by Hilton \cite{H} and is the first cross-intersection result of the kind we are dealing with.

%This brings us to the main result in this paper.

\subsection{The main result}

As pointed out above, ${[n] \choose r} = \mathcal{H}^{(r)}$ with $\mathcal{H} = 2^{[n]}$. %Recall that if $\mathcal{H} = 2^{[n]}$, then ${[n] \choose r} = \mathcal{H}^{(r)}$ and $\mu(\mathcal{H}) = n$. 
Thus, the theorem of Wang and Zhang %\cite{WZ2} %Theorem~\ref{WZthm} 
deals with the $r$th level and the $s$th level of the hereditary family $2^{[n]}$. %addresses the problem of determining $M(\mathcal{H}^{(r)},\mathcal{H}^{(s)},t)$ with $\mathcal{H}$ being the hereditary family $2^{[n]}$. 
We characterize the pairs in $M(\mathcal{H}^{(r)},\mathcal{H}^{(s)},t)$ for any hereditary family $\mathcal{H}$ with $\mu(\mathcal{H})$ sufficiently large depending on $r$, $s$, and $t$.
%For a family $\mathcal{F}$ and positive integers $r$, $s$, and $t$, let $C(\mathcal{F},r,s,t)$ denote the set $\{(\mathcal{A},\mathcal{B}) \colon \emptyset \neq \mathcal{A} \subseteq \mathcal{F}^{(r)}, \emptyset \neq \mathcal{B} \subseteq \mathcal{F}^{(s)},$ $\mathcal{A}$ and $\mathcal{B}$ are cross-$t$-intersecting, and if  $r=s$, then $|\mathcal{A}| \leq |\mathcal{B}|\}$, let $m(\mathcal{F},r,s,t) = \max\{|\mathcal{A}| + |\mathcal{B}| \colon (\mathcal{A},\mathcal{B}) \in C(\mathcal{F},r,s,t)\}$, and let $M(\mathcal{F},r,s,t)$ denote the set of pairs $(\mathcal{A}, \mathcal{B})$ in $C(\mathcal{F},r,s,t)$ such that $|\mathcal{A}| + |\mathcal{B}|$ is maximum, that is, $M(\mathcal{F},r,s,t) = \{(\mathcal{A},\mathcal{B}) \in C(\mathcal{F},r,s,t) \colon  |\mathcal{A}| + |\mathcal{B}| = m(\mathcal{F},r,s,t)\}$.

The paper \cite{Borg9} features the following two conjectures for $t=1$. %We now conveniently take $r \leq s$ as in \cite{Borg9}.

\begin{conj}[Weak Form \cite{Borg9}] \label{conj1} If $1 \leq r \leq s$ and $\mathcal{H}$ is a hereditary family with $\mu(\mathcal{H}) \geq r+s$, then for some $(\mathcal{A}, \mathcal{B}) \in M(\mathcal{H}^{(r)},\mathcal{H}^{(s)},1)$, $\mathcal{A}$ is a trivial $1$-intersecting family.
\end{conj}

\begin{conj}[Strong Form \cite{Borg9}] \label{conj2} If $1 \leq r \leq s$ and $\mathcal{H}$ is a hereditary family with $\mu(\mathcal{H}) \geq r+s$, then there exists a set $I$ in $\mathcal{H}$ such that $1 \leq |I| \leq r$ and for some $(\mathcal{A}, \mathcal{B}) \in M(\mathcal{H}^{(r)},\mathcal{H}^{(s)},1)$, $\mathcal{A} = \mathcal{H}^{(r)}(I)$ and $\mathcal{B} = \{B \in \mathcal{H}^{(s)} \colon B \cap I \neq \emptyset\}$.
\end{conj}
%
%A result in \cite{HM} (which is perhaps the first cross-intersection %result of the kind we are considering) says that the above conjecture %is true for $\mathcal{H} = 2^{[n]}$
Generalizing the above-mentioned result of Frankl and Tokushige \cite{FT1}, the main result in \cite{Borg9} tells us that for certain hereditary families $\mathcal{H}$, Conjecture~\ref{conj2} holds with $|I| = r$, in which case $\mathcal{A}$ consists of $I$ only and $\mathcal{B}$ consists of all the sets in $\mathcal{H}^{(s)}$ intersecting $I$. A question that arises immediately is whether this holds for every hereditary family. 
%we can improve Conjecture~\ref{conj2} to having $|I| = r$ so that . %of $\mathcal{H}^{(r)}$ and $\mathcal{H}^{(s)}$, respectively. 
This is answered in the negative in \cite{Borg9} too; \cite[Proposition~2.1]{Borg9} tells us that for any $2 \leq r \leq s$ and $n \geq r+s$, there are hereditary families $\mathcal{H}$ such that $\mu(\mathcal{H}) = n$ and no $(\mathcal{A}, \mathcal{B})$ in $M(\mathcal{H}^{(r)},\mathcal{H}^{(s)},1)$ satisfies Conjecture~\ref{conj2} with $|I| = r$.

Throughout the paper, %for any positive integers $r$, $s$ and $t$, 
we take
\begin{equation} c(r,s,t) = r + (s-t)\max \left\{2{s \choose t}, \; 2^r(r-t){r \choose t} + 1 \right\}. \nonumber
\end{equation} 
 
Note that Conjecture~\ref{conj2} is significantly stronger than Conjecture~\ref{conj1}. %We prove the more general cross-$t$-intersection problem for $\mu(\mathcal{H}) \geq n_{\rm N}(r,s,t)$, hence verifying Conjecture~\ref{conj2} for $\mu(\mathcal{H}) \geq n_{\rm N}(r,s,1)$. 
In Section~\ref{nonemptysection}, we prove the following generalization for $M(\mathcal{H}^{(r)},\mathcal{H}^{(s)},t)$ with $\mu(\mathcal{H}) \geq c(r,s,t)$, hence verifying Conjecture~\ref{conj2} for $\mu(\mathcal{H}) \geq c(r,s,1)$. 
%We prove an even stronger statement for cross-$t$-intersecting families for $\mu(\mathcal{H}) \geq n_{\rm N}(r,s,t)$, hence verifying Conjecture~\ref{conj2} for $\mu(\mathcal{H}) \geq n_{\rm N}(r,s,1)$.

%\begin{theorem} \label{nonemptyfam} If $1 \leq t \leq r \leq s$, $\mathcal{H}$ is a hereditary family with $\mu(\mathcal{H}) \geq c(r,s,t)$, and $(\mathcal{A}, \mathcal{B}) \in M(\mathcal{H}^{(r)},\mathcal{H}^{(s)},t)$, then there exists a set $I$ in $\mathcal{H}$ such that $t \leq |I| \leq r$, $\mathcal{A} = \mathcal{H}^{(r)}(I)$, and $\mathcal{B} = \{B \in \mathcal{H}^{(s)} \colon |B \cap I| \geq t\}$.
% 
%\[\mbox{$\mathcal{A} = \{A \in \mathcal{H}^{(r)} \colon I \subseteq A\}$, and $\mathcal{B} = \{B \in \mathcal{H}^{(s)} \colon |B \cap I| \geq t\}.$}\]
%\end{theorem}
%
%This is equivalent to the following statement.

\begin{theorem} \label{nonemptyfam} If $1 \leq t \leq r \leq s$, $\mathcal{H}$ is a hereditary family with $\mu(\mathcal{H}) \geq c(r,s,t)$, and $(\mathcal{A}, \mathcal{B}) \in M(\mathcal{H}^{(r)},\mathcal{H}^{(s)},t)$, then for some set $I$ in $\mathcal{H}$ with $t \leq |I| \leq r$, either %exactly one of the following holds:
%
%\begin{enumerate}[label=(\roman{*})]
%\item $\mathcal{A} = \mathcal{H}^{(r)}(I)$ and $\mathcal{B} = \{B \in \mathcal{H}^{(s)} \colon |B \cap I| \geq t\}$,
%\item $t < |I|$, $r = s$, $\mathcal{A} = \mathcal{H}^{(r)}(I)$, and $\mathcal{B} = \{B \in \mathcal{H}^{(s)} \colon |B \cap I| \geq t\}$.
%\end{enumerate}
\[\mbox{$\mathcal{A} = \mathcal{H}^{(r)}(I)$ and $\mathcal{B} = \{B \in \mathcal{H}^{(s)} \colon |B \cap I| \geq t\}$,}\] 
or 
\[\mbox{$r = s$, $t < |I|$, $\mathcal{A} = \{A \in \mathcal{H}^{(r)} \colon |A \cap I| \geq t\}$, and $\mathcal{B} = \mathcal{H}^{(s)}(I)$.}\] 
\end{theorem}
It immediately follows that 
\begin{equation} (\mathcal{H}^{(r)}(I), \{B \in \mathcal{H}^{(s)} \colon |B \cap I| \geq t\}) \in M(\mathcal{H}^{(r)},\mathcal{H}^{(s)},t)
\end{equation} 
(with $I$ as in Theorem~\ref{nonemptyfam}). Thus, the following holds.

\begin{theorem} If $1 \leq t \leq r \leq s$ and $\mathcal{H}$ is a hereditary family with $\mu(\mathcal{H}) \geq c(r,s,t)$, then
\[m(\mathcal{H}^{(r)},\mathcal{H}^{(s)},t) = |\mathcal{H}^{(r)}(I)| + |\{B \in \mathcal{H}^{(s)} \colon |B \cap I| \geq t\}|\]
for some set $I$ in $\mathcal{H}$ with $t \leq |I| \leq r$.
\end{theorem}

%\begin{prop} Theorem~\ref{nonemptyfam} and Theorem~\ref{nonemptyfamequiv} are equivalent.
%\end{prop}
%
%\textbf{Proof.} Clearly, Theorem~\ref{nonemptyfamequiv} implies Theorem~\ref{nonemptyfam}. Assume Theorem~\ref{nonemptyfam} and let $\mathcal{A}$ and $\mathcal{B}$ be as in Theorem~\ref{nonemptyfamequiv}. 

\begin{problem} For $1 \leq t \leq r \leq s$, let $\eta(r,s,t)$ be the smallest integer $n$ such that for every hereditary family $\mathcal{H}$ with $\mu(\mathcal{H}) \geq n$, $(\mathcal{H}^{(r)}(I), \{B \in \mathcal{H}^{(s)} \colon |B \cap I| \geq t\}) \in M(\mathcal{H}^{(r)}, \mathcal{H}^{(s)},t)$ for some $I \in \mathcal{H}$ with $t \leq |I| \leq r$. What is the value of $\eta(r,s,t)$?
\end{problem} 
By Theorem~\ref{nonemptyfam}, $\eta(r,s,t) \leq c(r,s,t)$. Clearly, for $\mathcal{H} = 2^{[n]}$, we have $\mu(\mathcal{H}) = n$, and $\mathcal{H}^{(r)}$ and $\mathcal{H}^{(s)}$ are cross-$t$-intersecting if and only if $n \leq r + s - t$. Thus, $\eta(r,s,t) \geq r + s - t + 1$. We conjecture that equality holds.

\begin{conj} \label{nonemptyconj} For $1 \leq t \leq r \leq s$, $\eta(r,s,t) = r + s - t + 1$.
\end{conj}

A \emph{graph} $G$ is a pair $(V,\mathcal{E})$ with $\mathcal{E} \subseteq {V \choose 2}$, and a subset $S$ of $V$ is called an \emph{independent set of $G$} if $\{i,j\} \notin \mathcal{E}$ for every $i, j \in S$. Let $\mathcal{I}_G$ denote the family of all independent sets of a graph $G$. The EKR problem for $\mathcal{I}_G$ was introduced in \cite{HT} and inspired many results \cite{BH1,BH2,HST,HT,HK,Wr}. %HS,  %There are many instances in which an %EKR-type result was obtained for a family $\mathcal{F}$ that is %actually the family $\mathcal{I}_G$ 
Many EKR-type results can be phrased in terms of independent sets of graphs; see \cite[page 2878]{BH2}. Clearly, $\mathcal{I}_G$ is a hereditary family. Kamat \cite{Kamat} %made the following conjecture. 
conjectured that if $\mu(\mathcal{I}_G) \geq 2r$, and $\mathcal{A}$ and $\mathcal{B}$ are cross-intersecting subfamilies of ${\mathcal{I}_G}^{(r)}$, then $|\mathcal{A}| + |\mathcal{B}| \leq |{\mathcal{I}_G}^{(r)}|$. We suggest the following strong generalization.

%\begin{conj}[\cite{Kamat}] If $G$ is a graph, $\mu(\mathcal{I}_G) \geq 2r$, and $\mathcal{A}$ and $\mathcal{B}$ are cross-intersecting subfamilies of ${\mathcal{I}_G}^{(r)}$, then $|\mathcal{A}| + |\mathcal{B}| \leq |{\mathcal{I}_G}^{(r)}|$.
%\end{conj}
%

\begin{conj}\label{nonemptyconjcor} If $1 \leq t \leq r \leq s$, $\mathcal{H}$ is a hereditary family with $\mu(\mathcal{H}) \geq r+s-t+1$, $\mathcal{A} \subseteq \mathcal{H}^{(r)}$, $\mathcal{B} \subseteq \mathcal{H}^{(s)}$, and $\mathcal{A}$ and $\mathcal{B}$ are cross-$t$-intersecting, then $|\mathcal{A}| + |\mathcal{B}| \leq |{\mathcal{H}}^{(s)}|$.
\end{conj}
%\begin{conj}\label{kamatgen} If $\mathcal{H}$ is a hereditary family, $\mu(\mathcal{H}) \geq 2r$, and $\mathcal{A}$ and $\mathcal{B}$ are cross-intersecting subfamilies of $\mathcal{H}^{(r)}$, then $|\mathcal{A}| + |\mathcal{B}| \leq |\mathcal{H}^{(r)}|$.
%\end{conj}
%
%Naturally, similarly to Conjecture~\ref{AK gen}, we also conjecture that if $\mu(\mathcal{H}) > 2r$, then the bound is attained if and only if one of $\mathcal{A}$ and $\mathcal{B}$ is $\mathcal{H}^{(r)}$ and the other is empty. 
%This is true for $\mu(\mathcal{H})$ sufficiently large depending only on $r$.
%
In other words, we conjecture that for $\mu(\mathcal{H}) \geq r+s-t+1$, if the cross-$t$-intersecting families $\mathcal{A}$ and $\mathcal{B}$ are allowed to be empty, then their sum of sizes is maximum if $\mathcal{A}$ is empty and $\mathcal{B}$ is $\mathcal{H}^{(s)}$. 

In Section~\ref{propertysection}, we establish some key properties of hereditary families that enable us to prove Theorem~\ref{nonemptyfam} and the following result.

\begin{lemma}\label{nonemptylemma} If $1 \leq t \leq r \leq s$, $\mathcal{H}$ is a hereditary family with $\mu(\mathcal{H}) \geq r+s-t+1$, $I$ is a set in $\mathcal{H}$ with $t \leq |I| \leq r$, $\mathcal{A} = \mathcal{H}^{(r)}(I)$, and $\mathcal{B} = \{B \in \mathcal{H}^{(s)} \colon |B \cap I| \geq t\}$, then $|\mathcal{A}| + |\mathcal{B}| \leq |\mathcal{H}^{(s)}|$, 
%
%\[|\mathcal{A}| + |\mathcal{B}| \leq |\mathcal{H}^{(s)}|,\]
%
%Moreover, if $\mu(\mathcal{H}) > r+s-t+1$ and $|\mathcal{A}| + |\mathcal{B}| = |\mathcal{H}^{(s)}|$, then $(\mathcal{A},\mathcal{B}) = (\emptyset, \mathcal{H}^{(s)})$.
and equality holds only if $t= 1$ and $\mu(\mathcal{H}) = r+s$.
\end{lemma}
Lemma~\ref{nonemptylemma} is also proved in Section~\ref{propertysection}. %It tells us that for $\mu(\mathcal{H}) \geq r+s-t+1$, if the cross-$t$-intersecting families are allowed to be empty, then the cross-$t$-intersecting families $\{A \in \mathcal{H}^{(r)} \colon I \subseteq A\}$ and $\{B \in \mathcal{H}^{(s)} \colon |B \cap I| \geq t\}$ (with $t \leq |I| \leq r$) %configuration in Theorem~\ref{nonemptyfam} 
%do not have a larger sum of sizes than the cross-$t$-intersecting families $\emptyset$ and $\mathcal{H}^{(s)}$. 
%Consequently, we have the following two results.
It immediately gives us the following.
\begin{theorem} If Conjecture~\ref{nonemptyconj} is true, then Conjecture~\ref{nonemptyconjcor} is true.
\end{theorem}
Together with Theorem~\ref{nonemptyfam}, Lemma~\ref{nonemptylemma} also immediately yields the following.

\begin{theorem} \label{nonemptyfamcor} If $1 \leq t \leq r \leq s$, $\mathcal{H}$ is a hereditary family with $\mu(\mathcal{H}) \geq c(r,s,t)$, $\mathcal{A} \subseteq \mathcal{H}^{(r)}$, $\mathcal{B} \subseteq \mathcal{H}^{(s)}$, and $\mathcal{A}$ and $\mathcal{B}$ are cross-$t$-intersecting, then $(\mathcal{A},\mathcal{B}) = (\emptyset,\mathcal{H}^{(s)})$ or $r = s$ and $(\mathcal{A},\mathcal{B}) = (\mathcal{H}^{(r)},\emptyset)$.
\end{theorem}
Therefore, Conjecture~\ref{nonemptyconjcor} is true if $\mu(\mathcal{H}) \geq c(r,s,t)$, and hence Kamat's conjecture is true if $\mu(\mathcal{I}_G) \geq c(r,r,1)$.

We mention that the analogous problem for cross-intersecting subfamilies of $\mathcal{H}$ is solved in \cite{Borg5}.

We now start working towards proving Theorem~\ref{nonemptyfam} and Lemma~\ref{nonemptylemma}. %The rest of the paper is dedicated to the proof of this result. %, concluded in Section~\ref{nonemptysection}. The tools are given in the next two sections.

\section{Key properties of hereditary families} \label{propertysection}

Hereditary families exhibit undesirable phenomena; see, for example, \cite[Example~1]{Borg}. The complete absence of symmetry makes intersection problems like the ones described above very difficult to deal with. Many of the well-known techniques in extremal set theory, such as the shifting technique (see \cite{F}), fail to work for hereditary families. The lemmas in this section and the next are the tools that will enable us to overcome such difficulties.
%The ingredients that will enable us to overcome such difficulties are given in this section and the next.

The two results below establish the properties of hereditary families that are fundamental to our work. The first one is given by {\cite[Corollary~3.2]{Borg}}. %established in \cite{Borg}. %We will need it for all our results.

\begin{lemma}[{\cite{Borg}}]\label{Spernercor} If $\mathcal{H}$ is a hereditary family and $0 \leq r \leq s \leq \mu(\mathcal{H}) - r$, then
\[|\mathcal{H}^{(s)}| \geq \frac{{\mu(\mathcal{H}) - r \choose
s - r}}{{s \choose s-r}}|\mathcal{H}^{(r)}|.\]
\end{lemma}

\begin{lemma}\label{mulemma} If $\mathcal{H}$ is a hereditary family, $X \subseteq Y$, $\mathcal{G}$ is the family $\{H \in \mathcal{H} \colon H \cap Y = X\}$, and $\mathcal{G} \neq \emptyset$, then \[\mu(\{G \backslash X \colon G \in \mathcal{G}\}) \geq \mu(\mathcal{H}) - |Y|.\]
\end{lemma}
\textbf{Proof.} Let $\mathcal{F} = \{G \backslash X \colon G \in \mathcal{G}\}$. Since $\mathcal{G} \neq \emptyset$, $\mathcal{F} \neq \emptyset$. Let $B$ be a base of $\mathcal{F}$ of size $\mu(\mathcal{F})$. Let $C = B \cup X$. Then $C \in \mathcal{G}$, and hence $C \in \mathcal{H}$. Let $D$ be a base of $\mathcal{H}$ such that $C \subseteq D$. Then $X \subseteq D$. 
%
%Suppose that $|D| \geq |B| + |Y| + 1$. 
Let $E = (D \backslash Y) \cup X$. Since $\mathcal{H}$ is hereditary and $E \subseteq D \in \mathcal{H}$, $E \in \mathcal{H}$. Since $E \cap Y = X$, $E \in \mathcal{G}$. Let $F = E \backslash X$. Then $F \in \mathcal{F}$. Since $C \subseteq D$ and $C \cap Y = E \cap Y = X$, $B \subseteq F$. Since $B$ is a base of $\mathcal{F}$, $B = F$. Thus, we have $\mu(\mathcal{F}) = |B| = |F| = |E| - |X| = |D \backslash Y| \geq |D| - |Y| \geq \mu(\mathcal{H}) - |Y|$.~\hfill{$\Box$} \\

For $X = Y$, the lemma above holds even if the family is not hereditary.

\begin{lemma} \label{mucor} If $\mathcal{F}$ is a family and $X$ is a set such that $\mathcal{F}(X) \neq \emptyset$, then
\[\mu(\{F \backslash X \colon F \in \mathcal{F} ( X )\}) \geq \mu(\mathcal{F}) - |X|.\]
\end{lemma}
\textbf{Proof.} Let $\mathcal{G} = \{F \backslash X \colon F \in \mathcal{F} ( X )\}$. Let $B$ be a base of $\mathcal{G}$ of size $\mu(\mathcal{G})$. Then $B \cup X$ is a base of $\mathcal{F}$. Thus, $\mu(\mathcal{F}) \leq |B| + |X| = \mu(\mathcal{G}) + |X|$.~\hfill{$\Box$}

\begin{lemma}\label{mainlemma2} If $0 \leq t \leq u \leq r$, $s \geq r + t - u$, $\mathcal{H}$ is a hereditary family with $\mu(\mathcal{H}) \geq r + s - t$, and $T$ is a $t$-element subset of a $u$-element set $U$ such that $\mathcal{H}^{(r)}(U) \neq \emptyset$, then 
\[|\{H \in \mathcal{H}^{(s)} \colon H \cap U = T\}| \geq \frac{ {\mu(\mathcal{H}) - r \choose s + u - r - t} }{ {s - t \choose s + u - r - t} } |\mathcal{H}^{(r)}(U)|.\]
\end{lemma}
\textbf{Proof.} Let $\mathcal{S} = \{H \in \mathcal{H}^{(s)} \colon H \cap U = T\}$. Since $\mathcal{H}^{(r)}(U) \neq \emptyset$, $\mathcal{H}(U) \neq \emptyset$. Let $\mathcal{I} = \{H \backslash U \colon H \in \mathcal{H}(U)\}$. Since $\mathcal{H}$ is hereditary, $\mathcal{I}$ is hereditary. By Lemma~\ref{mucor}, $\mu(\mathcal{I}) \geq \mu(\mathcal{H}) - u$. Let $p = r - u$ and $q = s - t$. Since $\mu(\mathcal{H}) \geq r + s - t$, $\mu(\mathcal{I}) \geq r + s - t - u = p + q$. We have $0 \leq p \leq q \leq \mu(\mathcal{I}) - p$. Therefore, by Lemma~\ref{Spernercor},
\begin{gather} |\mathcal{I}^{(q)}| \geq
\frac{ {\mu(\mathcal{I})-p \choose q - p} }{ {q \choose q-p} } |\mathcal{I}^{(p)}|. \label{mainlemma2.1}
\end{gather}
Clearly, $|\mathcal{I}^{(p)}| = |\mathcal{H}^{(r)}(U)|$. Consider any $A \in \mathcal{I}^{(q)}$. Since $A \cup T \subseteq A \cup U \in \mathcal{H}(U)$ and $\mathcal{H}$ is hereditary, $A \cup T \in \mathcal{H}$. Since $|A \cup T| = s$ and $(A \cup T) \cap U = T$, it follows that $A \cup T \in \mathcal{S}$. Thus, $|\mathcal{I}^{(q)}| \leq |\mathcal{S}|$. Therefore, by (\ref{mainlemma2.1}),
\begin{align} |\mathcal{S}| &\geq \frac{ {\mu(\mathcal{I})-p \choose q - p} }{ {q \choose q - p} } |\mathcal{H}^{(r)}(U)| \geq \frac{ {(\mu(\mathcal{H}) - u) - (r - u) \choose (s-t) - (r-u)} }{ {s - t \choose (s-t) - (r-u)} } |\mathcal{H}^{(r)}(U)| = \frac{ {\mu(\mathcal{H}) - r \choose s + u - r - t} }{ {s - t \choose s + u - r - t} } |\mathcal{H}^{(r)}(U)|, \nonumber
\end{align}
as required.~\hfill{$\Box$}\\
\\
\textbf{Proof of Lemma~\ref{nonemptylemma}.} Let $t' = t-1$. For each $T \in {I \choose t'}$, let $\mathcal{S}_T = \{H \in \mathcal{H}^{(s)} \colon H \cap I = T\}$. Consider any $T \in {I \choose t'}$. We have $\mathcal{S}_T \cap \mathcal{B} = \emptyset$. Also, by Lemma~\ref{mainlemma2}, 
\[|\mathcal{S}_T| \geq \frac{ {\mu(\mathcal{H}) - r \choose s + |I| - r - t'} }{ {s - t' \choose s + |I| - r - t'} } |\mathcal{H}^{(r)}(I)| \geq \frac{ {s-t+1 \choose s + |I| - r - t + 1} }{ {s - t + 1 \choose s + |I| - r - t + 1}} |\mathcal{H}^{(r)}(I)| = |\mathcal{A}|,\]
and equality holds throughout only if $\mu(\mathcal{H}) = r+s-t+1$. We have $|\mathcal{H}^{(s)}| \geq |\mathcal{B} \cup \bigcup_{T \in {I \choose t'}} \mathcal{S}_T| = |\mathcal{B}| + \sum_{T \in {I \choose t'}} |\mathcal{S}_T| \geq |\mathcal{B}| + {|I| \choose t'}|\mathcal{A}|\geq |\mathcal{A}| + |\mathcal{B}|$, and equality holds throughout only if $\mu(\mathcal{H}) = r+s-t+1$ and $t' = 0$. The result follows.~\hfill{$\Box$}

\section{Proof of Theorem~\ref{nonemptyfam}}\label{nonemptysection}

If a set $X$ $t$-intersects each set in a family $\mathcal{A}$, then we call $X$ a \emph{$t$-transversal of $\mathcal{A}$}. 

\begin{lemma}\label{transversalemma1} If $X$ is a $t$-transversal of a family $\mathcal{A}$, then
\[|\mathcal{A}| \leq {|X| \choose t} |\mathcal{A}(T)|\]
for some $T \in {X \choose t}$.
\end{lemma}
\textbf{Proof.} Let $\mathcal{X} = {X \choose t}$. Let $T \in {X \choose t}$ such that $|\mathcal{A}(I)| \leq |\mathcal{A}(T)|$ for each $I \in \mathcal{X}$. Since $|A \cap X| \geq t$ for each $A \in \mathcal{A}$, we clearly have $\mathcal{A} = \bigcup_{I \in \mathcal{X}} \mathcal{A}(I)$. Thus, $|\mathcal{A}| = \left| \bigcup_{I \in \mathcal{X}}\mathcal{A}(I) \right| \leq \sum_{I \in \mathcal{X}} |\mathcal{A}(I)| \leq \sum_{I \in \mathcal{X}}|\mathcal{A}(T)| = |\mathcal{X}| |\mathcal{A}(T)| = {|X| \choose t} |\mathcal{A}(T)|$.~\hfill{$\Box$}

\begin{lemma} \label{transversalemma2} If $X$ is a $t$-transversal of a family $\mathcal{A}$, $T$ is a set of size $t$, and $T \nsubseteq X$, then 
\[\mathcal{A}(T) = \bigcup_{x \in X \backslash T} \mathcal{A}(T \cup \{x\}).\]
\end{lemma}
\textbf{Proof.} Obviously, $\bigcup_{x \in X \backslash T} \mathcal{A}(T \cup \{x\}) \subseteq \mathcal{A}(T)$. For each $A \in \mathcal{A}$, we have 
\[t \leq |A \cap X| = |A \cap (X \cap T)| + |A \cap (X \backslash T)| \leq t-1 + |A \cap (X \backslash T)|\] 
(as $|T| = t$ and $T \nsubseteq X$), and hence $|A \cap (X \backslash T)| \geq 1$. Thus, for each $A \in \mathcal{A}(T)$, we have $a \in A$ for some $a \in X \backslash T$, and hence $A \in \mathcal{A}(T \cup \{a\}) \subseteq \bigcup_{x \in X \backslash T} \mathcal{A}(T \cup \{x\})$. Therefore, we have $\mathcal{A}(T) \subseteq \bigcup_{x \in X \backslash T} \mathcal{A}(T \cup \{x\}) \subseteq \mathcal{A}(T)$. The result follows.~\hfill{$\Box$}

\begin{lemma} \label{transversalemma3} If $\mathcal{A}$ and $\mathcal{B}$ are non-empty cross-$t$-intersecting families such that $\mathcal{A}$ is $r$-uniform, $\mathcal{B}$ is $s$-uniform, and $\mathcal{B}$ is not a trivial $t$-intersecting family, then there exist $B, X \in \mathcal{B}$ such that
\[|\mathcal{A}| \leq s{s \choose t} |\mathcal{A}(T \cup \{x\})|\]
for some $T \in {B \choose t}$ and some $x \in X \backslash T$.
\end{lemma}
\textbf{Proof.} Since $\mathcal{A}$ and $\mathcal{B}$ are cross-$t$-intersecting, each set in $\mathcal{A}$ is a $t$-transversal of $\mathcal{B}$, and each set in $\mathcal{B}$ is a $t$-transversal of $\mathcal{A}$. Let $B \in \mathcal{B}$. By Lemma~\ref{transversalemma1}, $|\mathcal{A}| \leq {|B| \choose t} |\mathcal{A}(T)| = {s \choose t}|\mathcal{A}(T)|$ for some $T \in {B \choose t}$. Since $\mathcal{B}$ is not a trivial $t$-intersecting family, $T \nsubseteq X$ for some $X \in \mathcal{B}$. By Lemma~\ref{transversalemma2}, $\mathcal{A}(T) = \bigcup_{x \in X \backslash T} \mathcal{A}(T \cup \{x\})$, so $|\mathcal{A}(T)| \leq \sum_{x \in X \backslash T} |\mathcal{A}(T \cup \{x\})|$. Let $x^* \in X \backslash T$ such that $|\mathcal{A}(T \cup \{x\})| \leq |\mathcal{A}(T \cup \{x^*\})|$ for each $x \in X \backslash T$. Let $Y = T \cup \{x^*\}$. Thus, $|\mathcal{A}(T)| \leq \sum_{x \in X \backslash T} |\mathcal{A}(Y)| = |X \backslash T| |\mathcal{A}(Y)| \leq s |\mathcal{A}(Y)|$, and hence $|\mathcal{A}| \leq {s \choose t}s|\mathcal{A}(Y)|$.~\hfill{$\Box$}
\begin{lemma} \label{transversalemma4} If $1 \leq t \leq r$, $\mathcal{H}$ is a hereditary family with $\mu(\mathcal{H}) \geq 2r-t$, $\emptyset \neq \mathcal{A} \subseteq \mathcal{H}^{(r)}$, $\mathcal{B}$ is a non-empty $s$-uniform family that is not a trivial $t$-intersecting family, and $\mathcal{A}$ and $\mathcal{B}$ are cross-$t$-intersecting, then there exists a $t$-element set $T$ such that
\[|\mathcal{A}| < \frac{s(r-t)}{\mu(\mathcal{H})-r}{s \choose t} |\mathcal{H}^{(r)}(T)|\]
and $T \subseteq B$ for some $B \in \mathcal{B}$.
\end{lemma}
\textbf{Proof.} By Lemma~\ref{transversalemma3}, there exist $B, X \in \mathcal{B}$ such that such that $|\mathcal{A}| \leq s{s \choose t} |\mathcal{A}(T \cup \{x\})|$ for some $T \in {B \choose t}$ and some $x \in X \backslash T$.  Since $\mathcal{A} \neq \emptyset$, it follows that $\mathcal{A}(T \cup \{x\}) \neq \emptyset$, so $\mathcal{H}^{(r)}(T \cup \{x\}) \neq \emptyset$. Let $\mathcal{G} = \{H \in \mathcal{H}^{(r)} \colon H \cap (T \cup \{x\}) = T\}$. We have $|\mathcal{A}(T \cup \{x\})| \leq |\mathcal{H}^{(r)}(T \cup \{x\})| \leq \frac{r-t}{\mu(\mathcal{H})-r} |\mathcal{G}|$ by Lemma~\ref{mainlemma2}. Since $|\mathcal{H}^{(r)}(T)| = |\mathcal{G}| + |\mathcal{H}^{(r)}(T \cup \{x\})| > |\mathcal{G}|$, we obtain $|\mathcal{A}(T \cup \{x\})| < \frac{r-t}{\mu(\mathcal{H})-r} |\mathcal{H}^{(r)}(T)|$. Since $|\mathcal{A}| \leq s{s \choose t} |\mathcal{A}(T \cup \{x\})|$, the result follows.~\hfill{$\Box$} \\

We now settle a few calculations so that in the formal proof of the theorem we can focus on the combinatorial argument.

\begin{prop}\label{calc} If $1 \leq t \leq r \leq s$, $(r,s) \neq (t,t)$, and $n \geq c(r,s,t)$, then
\begin{flalign*} \mbox{(i) } \; &\frac{r(s-t)}{n - s} {r \choose t} < \frac{1}{2}. &\nonumber \\
\mbox{(ii) } \; &{s \choose t} \leq \frac{1}{2} \frac{{n-r \choose s-r}}{{s-t \choose s-r}} \mbox{ if } r < s. &\nonumber 
\end{flalign*}
\end{prop}
\textbf{Proof.} By straightforward induction, $2^a \geq 2a$ for every positive integer $a$. Since $t \leq r \leq s$ and $(r,s) \neq (t,t)$, either $t < r$ or $t = r < s$. If $t < r$, then, since $n \geq 2^r(r-t)(s-t){r \choose t} + r + s - t$, we have $n > 2r(s-t){r \choose t} + s$, which yields (i). If $t = r < s$, then, since $n \geq 2(s-t){s \choose t} + r$, we have $n \geq 2(s-t){t+1 \choose t} + t = 2(t+1)(s-t) + t > 2t(s-t) + s = 2r(s-t){r \choose t} + s$ (as $r = t$), which yields (i). 

Suppose $s > r$. Then $s > t$. Since $n \geq 2(s-t){s \choose t} + r$, we have $n - r > s-t > 0$ and ${s \choose t} \leq \frac{1}{2} \left( \frac{n-r}{s-t} \right)$. Thus, ${s \choose t} \leq \frac{1}{2} \prod_{i=0}^{s-r-1} \left( \frac{n-r-i}{s-t-i} \right) = \frac{1}{2} \frac{{n-r \choose s-r}}{{s-t \choose s-r}}$, which confirms (ii).~\hfill{$\Box$}
\\
\\
\textbf{Proof of Theorem~\ref{nonemptyfam}.} Let $n = c(r,s,t)$. Let $\mathcal{A}$ and $\mathcal{B}$ be as in the theorem.\medskip %Since $\mathcal{A}$ and $\mathcal{B}$ are cross-$t$-intersecting, each set in $\mathcal{A}$ is a $t$-transversal of $\mathcal{B}$, and each set in $\mathcal{B}$ is a $t$-transversal of $\mathcal{A}$.\medskip %We investigate two complementary cases separately.
%We first show that $\mathcal{A}$ is a trivial $t$-intersecting family.\medskip

\textit{Case 1: $\mathcal{A}$ is a trivial $t$-intersecting family.} Let $I = \bigcap_{A \in \mathcal{A}} A$, $\mathcal{C} = \mathcal{H}^{(r)}(I)$, and $\mathcal{D} = \{H \in \mathcal{H}^{(s)} \colon |H \cap I| \geq t\}$. Then $t \leq |I| \leq r$, $I \in \mathcal{H}$ (as $\mathcal{H}$ is hereditary), and $\mathcal{A} \subseteq \mathcal{C}$.

Suppose $|I| = r$. Then $\mathcal{A} = \{I\}$ and, since $\mathcal{A}$ and $\mathcal{B}$ are cross-$t$-intersecting, $\mathcal{B} \subseteq \mathcal{D}$. Since $\{I\}$ and $\mathcal{D}$ are cross-$t$-intersecting, and since $(\mathcal{A}, \mathcal{B}) \in M(\mathcal{H}^{(r)},\mathcal{H}^{(s)},t)$, we obtain $\mathcal{B} = \mathcal{D}$, as required.

Now suppose $|I| < r$. Let $\mathcal{A}' = \{A \backslash I \colon A \in \mathcal{A}\}$, $\mathcal{I} = \{H \backslash I \colon H \in \mathcal{H}(I)\}$, and $r' = r-|I|$. Then $\mathcal{A}' \subseteq \mathcal{I}^{(r')}$, $\mathcal{I}$ is hereditary, and, by Lemma~\ref{mucor}, $\mu(\mathcal{I}) \geq \mu(\mathcal{H}) - |I|$. By the definition of $I$, $\bigcap_{E \in \mathcal{A}'} E = \emptyset$. Thus, $\mathcal{A}'$ is not a trivial $1$-intersecting family. For each $i \in \{0\} \cup [t-1]$, let $\mathcal{B}_i = \{B \in \mathcal{B} \colon |B \cap I| = i\}$. Let $\mathcal{B}_t = \{B \in \mathcal{B} \colon |B \cap I| \geq t\}$. Then $\mathcal{B} = \bigcup_{i=0}^t \mathcal{B}_i$. Let $J = \{i \in \{0\} \cup [t-1] \colon \mathcal{B}_i \neq \emptyset\}$. 

Suppose $J = \emptyset$. Then $\mathcal{B} = \mathcal{B}_t$. Hence $\mathcal{B} \subseteq \mathcal{D}$. Thus, as required, we obtain $\mathcal{A} = \mathcal{C}$ and $\mathcal{B} = \mathcal{D}$, because $\mathcal{A} \subseteq \mathcal{C}$, $\mathcal{C}$ and $\mathcal{D}$ are cross-$t$-intersecting, and $(\mathcal{A}, \mathcal{B}) \in M(\mathcal{H}^{(r)},\mathcal{H}^{(s)},t)$.

We now show that indeed $J = \emptyset$.

Suppose $J \neq \emptyset$. Consider any $j \in J$. For any $S \in {I \choose j}$, let $\mathcal{B}_{j,S} = \{B \in \mathcal{B}_j \colon B \cap I = S\}$. Then $\mathcal{B}_j = \bigcup_{S \in {I \choose j}} \mathcal{B}_{j,S}$. Let $\mathcal{S}_j = \{S \in {I \choose j} \colon \mathcal{B}_{j,S} \neq \emptyset\}$. Since $\mathcal{B}_j \neq \emptyset$, $\mathcal{S}_j \neq \emptyset$. Consider any $S \in \mathcal{S}_j$. Let $\mathcal{B}_{j,S}' = \{B \backslash S \colon B \in \mathcal{B}_{j,S}\}$, $\mathcal{H}_{j,S} = \{H \in \mathcal{H} \colon H \cap I = S\}$, $\mathcal{J}_{j,S} = \{H \backslash S \colon H \in \mathcal{H}_{j,S}\}$, $s_j = s - j$, and $t_j = t - j$. Then $\emptyset \neq \mathcal{B}_{j,S}' \subseteq {\mathcal{J}_{j,S}}^{(s_j)}$, $\mathcal{J}_{j,S}$ is hereditary, and, by Lemma~\ref{mulemma}, 
\[\mu(\mathcal{J}_{j,S}) \geq \mu(\mathcal{H}) - |I| > n - r \geq 2(s-t){s \choose t} \geq 2s(s-t) \geq 2s > 2s_j - t_j\] 
(note that $s > t$ as $t \leq |I| < r \leq s$). Since $\mathcal{A}$ and $\mathcal{B}$ are cross-$t$-intersecting, $\mathcal{A}'$ and $\mathcal{B}_{j,S}'$ are cross-$t_j$-intersecting. Since $t_j \geq 1$ and $\mathcal{A}'$ is not a trivial $1$-intersecting family, $\mathcal{A}'$ is not a trivial $t_j$-intersecting family. By Lemma~\ref{transversalemma4}, there exists a $t_j$-element set $X_{j,S}$ such that
\[|\mathcal{B}_{j,S}'| < \frac{r'(s_j-t_j)}{\mu(\mathcal{J}_{j,S}) - s_j} {r' \choose t_j} |{\mathcal{J}_{j,S}}^{(s_j)}(X_{j,S}) |\]
and $X_{j,S} \subseteq E_{j,S}$ for some $E_{j,S} \in \mathcal{A}'$. We have $|\mathcal{B}_{j,S}'| = |\mathcal{B}_{j,S}|$. Let $T_{j,S} = S \cup X_{j,S}$. Then $|{\mathcal{J}_{j,S}}^{(s_j)}(X_{j,S})| = |{\mathcal{H}_{j,S}}^{(s)}(T_{j,S})|$. Thus, 
\begin{align} |\mathcal{B}_{j,S}| &< \frac{r'(s_j-t_j)}{\mu(\mathcal{J}_{j,S}) - s_j} {r' \choose t_j}|{\mathcal{H}_{j,S}}^{(s)}(T_{j,S})| \leq \frac{(r-|I|)(s-t)}{\mu(\mathcal{H}) - |I| + j - s}{r-|I| \choose t-j}|{\mathcal{H}_{j,S}}^{(s)}(T_{j,S})|. \nonumber 
\end{align}
Since $\mathcal{A}'$ and $\mathcal{B}_{j,S}'$ are cross-$t_j$-intersecting, we have $r' \geq t_j$, that is, $r - |I| \geq t - j$. Since $0 \leq j \leq t-1$, $t \leq |I| \leq r-1$, and $\mu(\mathcal{H}) \geq n$, we therefore have
\begin{align} |\mathcal{B}_{j,S}| &< \frac{(r-t)(s-t)}{n + t - r - s} {r-j \choose t-j} |{\mathcal{H}_{j,S}}^{(s)}(T_{j,S})| \leq \frac{1}{2^r} |{\mathcal{H}_{j,S}}^{(s)}(T_{j,S})| \nonumber 
\end{align}
as $n \geq (r-t)(s-t)2^r{r \choose t} + r + s - t \geq (r-t)(s-t)2^r{r-j \choose t-j} + r + s - t$.
 
Let $j^* \in J$ and $S^* \in \mathcal{S}_{j^*}$ such that for each $j \in J$, $|{\mathcal{H}_{j,S}}^{(s)} ( T_{j,S} )| \leq |{\mathcal{H}_{j^*,S^*}}^{(s)} ( T_{j^*,S^*} )|$ for each $S \in \mathcal{S}_j$. We have
\begin{align} |\mathcal{B}| &= |\mathcal{B}_t| + \sum_{j \in J} |\mathcal{B}_j| \leq |\mathcal{D}| + \sum_{j \in J} \sum_{S \in \mathcal{S}_j}|\mathcal{B}_{j,S}| < |\mathcal{D}| + \sum_{j \in J} \sum_{S \in \mathcal{S}_j} \frac{1}{2^r} |{\mathcal{H}_{j,S}}^{(s)} ( T_{j,S} ) | \nonumber \\
&\leq |\mathcal{D}| + \sum_{j \in J} \sum_{S \in \mathcal{S}_j} \frac{1}{2^r} |{\mathcal{H}_{j^*,S^*}}^{(s)} ( T_{j^*,S^*} ) | \leq |\mathcal{D}| + \frac{1}{2^r} |{\mathcal{H}_{j^*,S^*}}^{(s)} ( T_{j^*,S^*} ) |\sum_{j \in J} \sum_{S \in \mathcal{S}_j} 1 \nonumber
\end{align}
and $\sum_{j \in J} \sum_{S \in \mathcal{S}_j} 1 = \sum_{j \in J} |\mathcal{S}_j| < \sum_{j = 0}^{|I|} {|I| \choose j} = 2^{|I|} \leq 2^{r-1}$. Thus,
\begin{equation} |\mathcal{B}| < |\mathcal{D}| + \frac{1}{2}|{\mathcal{H}_{j^*,S^*}}^{(s)} ( T_{j^*,S^*} )|. \label{18}
\end{equation}

For convenience, let $j = j^*$ and $S = S^*$. Let $B' \in \mathcal{B}_{j,S}'$. Recall that $\mathcal{A}'$ and $\mathcal{B}_{j,S}'$ are cross-$t_j$-intersecting, so $B'$ is a $t_j$-transversal of $\mathcal{A}'$. By Lemma~\ref{transversalemma1}, $|\mathcal{A}'| \leq {|B'| \choose t_j}|\mathcal{A}'(X^*)|$ for some $X^* \in {B' \choose t_j}$. Thus, we have
\begin{align} 0 < |\mathcal{A}| &= |\mathcal{A}'| \leq {s-j \choose t-j}|\mathcal{I}^{(r')}(X^*)| \leq {s \choose t}|\mathcal{I}^{(r')}(X^*)|. \label{20}
\end{align}
Let $\mathcal{K} = \{E \backslash X^* \colon E \in \mathcal{I} ( X^* )\}$, $p = r' - |X^*|$, and $q = s_j - |X^*|$. We have $p = r - |I| - t_j = r - |I| - t + j \leq r-t-1$ and $q = s_j - t_j = s - t \geq r-t \geq p+1$. Since $\mathcal{I}$ is hereditary, $\mathcal{K}$ is hereditary. Since $|\mathcal{I}(X^*)| \geq |\mathcal{I}^{(r')}(X^*)|$, $|\mathcal{I}(X^*)| > 0$ by (\ref{20}). Thus, by Lemma~\ref{mucor}, $\mu(\mathcal{K}) \geq \mu(\mathcal{I}) - |X^*| \geq \mu(\mathcal{H}) - |I| - t_j \geq n - |I| - t_j > r+s - |I| - t_j = p + s > p + q$. By Lemma~\ref{Spernercor}, 
\begin{equation} |\mathcal{K}^{(q)}| \geq \frac{{\mu(\mathcal{K}) - p \choose q-p}}{{q \choose q-p}}|\mathcal{K}^{(p)}| = |\mathcal{K}^{(p)}|\prod_{i = 0}^{q-p-1} \frac{\mu(\mathcal{K}) - p - i}{q-i} \geq |\mathcal{K}^{(p)}|\left(\frac{\mu(\mathcal{K}) - p}{q}\right)^{q-p}. \nonumber
\end{equation}
Since $q - p \geq 1$ and $\mu(\mathcal{K}) \geq n - |I| - t_j = n - r + p \geq p + 2(s-t){s \choose t} = p + 2q{s \choose t}$, $|\mathcal{K}^{(q)}| \geq 2{s \choose t}|\mathcal{K}^{(p)}|$. Thus, since $|\mathcal{K}^{(p)}| = |\mathcal{I}^{(r')}(X^*)|$ and $|\mathcal{K}^{(q)}| = |\mathcal{I}^{(s_j)}(X^*)|$, 
\begin{equation} {s \choose t} |\mathcal{I}^{(r')}(X^*)| \leq \frac{1}{2}|\mathcal{I}^{(s_j)}(X^*)|. \label{22}
\end{equation} 
Let $\mathcal{L} = \mathcal{H}^{(|I| + s_j)}(I \cup X^*)$.  Then $\mathcal{I}^{(s_j)}(X^*) = \{H \backslash I \colon H \in \mathcal{L}\}$. Let $\mathcal{L}' = \{L \backslash (I \backslash S) \colon L \in \mathcal{L}\}$. Since $\mathcal{H}$ is hereditary, $\mathcal{L}' \subseteq \mathcal{H}$. For each $H \in \mathcal{L}'$, we have $|H| = s_j + |I| - (|I|-|S|) = s$, $H \cap I = S$, and $S \cup X^* \subseteq H$. Thus, $\mathcal{L}' \subseteq {\mathcal{H}_{j,S}}^{(s)}(S \cup X^*)$. Let $T_1 = S \cup X^*$. We have $|\mathcal{I}^{(s_j)}(X^*)| = |\mathcal{L}| = |\mathcal{L}'| \leq |{\mathcal{H}_{j,S}}^{(s)}(T_1)|$. Together with (\ref{20}) and (\ref{22}), this gives us 
\begin{equation} |\mathcal{A}| \leq \frac{1}{2}|{\mathcal{H}_{j,S}}^{(s)}(T_1)|. \label{24}
\end{equation}
Let $T_2 = T_{j,S}$. Let $\mathcal{E}$ be a member of $\{{\mathcal{H}_{j,S}}^{(s)}(T_1), {\mathcal{H}_{j,S}}^{(s)}(T_2)\}$ of maximum size. Recall that above we set $j = j^*$ and $S = S^*$. By (\ref{18}) and (\ref{24}),
\begin{equation} |\mathcal{A}| + |\mathcal{B}| < \frac{1}{2}|{\mathcal{H}_{j,S}}^{(s)}(T_1)| + |\mathcal{D}| + \frac{1}{2}|{\mathcal{H}_{j,S}}^{(s)}(T_2)| \leq |\mathcal{D}| + |\mathcal{E}|. \label{25}
\end{equation}
Let
\[X' = \left\{ \begin{array}{ll} X^* & \mbox{if $\mathcal{E} =
{\mathcal{H}_{j,S}}^{(s)}(T_1)$;}\\
X_{j,S} & \mbox{if $\mathcal{E} =
{\mathcal{H}_{j,S}}^{(s)}(T_2)$.}
\end{array} \right.\]
Let $F = I \cup X'$. Let $\mathcal{F} = \mathcal{H}^{(r)}(F)$ and $\mathcal{G} = \mathcal{D} \cup \mathcal{E}$. If $X' = X^*$, then, since $|\mathcal{F}| = |\mathcal{H}^{(|I|+r')}( I \cup X^* )| = |\mathcal{I}^{(r')}( X^* )|$, $|\mathcal{F}| > 0$ by (\ref{20}). If $X' = X_{j,S}$, then, since $X_{j,S} \subseteq E_{j,S} \in \mathcal{A}'$, we have $F \subseteq I \cup E_{j,S} \in \mathcal{A}$, and hence $I \cup E_{j,S} \in \mathcal{F}$. Therefore, $\mathcal{F} \neq \emptyset$. By (\ref{25}), $\mathcal{G} \neq \emptyset$. For each $G \in \mathcal{D}$, $|G \cap F| \geq |G \cap I| \geq t$. For some $i \in [2]$,  $\mathcal{E} = {\mathcal{H}_{j,S}}^{(s)}(T_i)$ and $T_i = S \cup X'$; thus, for each $G \in \mathcal{E}$, $|G \cap F| \geq |T_i \cap F| = |S| + |X'| = j + t_j = t$. For every $G \in \mathcal{G}$ and every $H \in \mathcal{F}$, $|G \cap H| \geq |G \cap F|$, so $|G \cap H| \geq t$. Thus, $\mathcal{F}$ and $\mathcal{G}$ are cross-$t$-intersecting. For each $H \in \mathcal{E}$, $|H \cap I| = |S| = j < t$. Thus, $\mathcal{D} \cap \mathcal{E} = \emptyset$, and hence $|\mathcal{G}| = |\mathcal{D}| + |\mathcal{E}|$. Bringing all the pieces together, we have that $\emptyset \neq \mathcal{F} \subseteq \mathcal{H}^{(r)}$, $\emptyset \neq \mathcal{G} \subseteq \mathcal{H}^{(s)}$, $\mathcal{F}$ and $\mathcal{G}$ are cross-$t$-intersecting, and, by (\ref{25}), 
\[|\mathcal{A}| + |\mathcal{B}| < |\mathcal{G}| < |\mathcal{F}| + |\mathcal{G}|,\]
contradicting $(\mathcal{A}, \mathcal{B}) \in M(\mathcal{H}^{(r)},\mathcal{H}^{(s)},t)$.\medskip

\textit{Case 2: $\mathcal{A}$ is not a trivial $t$-intersecting family.} If $t = s$, then $t=r=s$ and $n = r = 2s-t$. If $t < s$, then $n > 2s$. Thus, $\mu(\mathcal{H}) \geq 2s-t$. By Lemma~\ref{transversalemma4}, there exists a $t$-element set $T_{\mathcal{B}}$ such that
\begin{equation} |\mathcal{B}| < \frac{r(s-t)}{\mu(\mathcal{H}) - s} {r \choose t} |\mathcal{H}^{(s)}(T_{\mathcal{B}})|. \label{14}
\end{equation}

Suppose $r < s$. Let $D \in \mathcal{B}$. Since $\mathcal{A}$ and $\mathcal{B}$ are cross-$t$-intersecting, $D$ is a $t$-transversal of $\mathcal{A}$. By Lemma~\ref{transversalemma1},
\begin{equation} |\mathcal{A}| \leq {|D| \choose t}|\mathcal{A}(T_D)| \leq {s \choose t}|\mathcal{H}^{(r)} (T_D)| \label{15}
\end{equation}
for some $T_D \in {D \choose t}$. Let $\mathcal{G} = \{H \backslash T_D \colon H \in \mathcal{H}( T_D )\}$. Then $\mathcal{G}$ is hereditary. Since $0 < |\mathcal{A}| \leq {s \choose t}|\mathcal{H}^{(r)}(T_D)| \leq {s \choose t}|\mathcal{H}(T_D)| = {s \choose t}|\mathcal{G}|$, $\mathcal{G} \neq \emptyset$. Thus, by Lemma~\ref{mucor}, $\mu(\mathcal{G}) \geq \mu(\mathcal{H}) - |T_D| = \mu(\mathcal{H}) - t$. By Lemma~\ref{Spernercor}, 
\[|\mathcal{G}^{(s-t)}| \geq \frac{{\mu(\mathcal{G}) - (r-t) \choose
(s-t) - (r-t)}}{{s-t \choose (s-t) - (r-t)}}|\mathcal{G}^{(r-t)}| = \frac{{\mu(\mathcal{G}) + t - r \choose
s-r}}{{s-t \choose s-r}}|\mathcal{G}^{(r-t)}|.\]
Clearly, $|\mathcal{H}^{(r)} ( T_{D} )| = |\mathcal{G}^{(r-t)}|$ and $|\mathcal{H}^{(s)} ( T_{D} )| = |\mathcal{G}^{(s-t)}|$. Let $T' \in \mathcal{H}^{(t)}$ such that $|\mathcal{H}^{(s)}( T )| \leq |\mathcal{H}^{(s)}( T' )|$ for all $T \in \mathcal{H}^{(t)}$. Since $\mathcal{A} \neq \emptyset$, $|\mathcal{H}^{(r)} ( T_{D} )| > 0$ by (\ref{15}). Since $\mathcal{H}$ is hereditary and $T_D$ is a $t$-element subset of every member of $\mathcal{H}^{(r)}(T_{D})$, we have $T_D \in \mathcal{H}^{(t)}$, and hence $|\mathcal{H}^{(s)}(T_{D})| \leq |\mathcal{H}^{(s)}(T')|$. Thus, we have
\begin{align} 0 < \frac{{\mu(\mathcal{H}) - r \choose
s-r}}{{s-t \choose s-r}}|\mathcal{H}^{(r)} ( T_{D} )| &\leq \frac{{\mu(\mathcal{G}) + t - r \choose
s-r}}{{s-t \choose s-r}}|\mathcal{H}^{(r)} ( T_{D} )| = \frac{{\mu(\mathcal{G}) + t - r \choose
s-r}}{{s-t \choose s-r}}|\mathcal{G}^{(r-t)}| \nonumber \\
&\leq |\mathcal{G}^{(s-t)}| = |\mathcal{H}^{(s)} ( T_{D} )| \leq |\mathcal{H}^{(s)}(T')|. \label{16a}
\end{align}
Thus, $\mathcal{H}^{(s)}(T') \neq \emptyset$. Since $\mathcal{H}$ is hereditary and every set in $\mathcal{H}^{(s)}(T')$ has an $r$-element subset containing $T'$, $\mathcal{H}^{(r)}(T') \neq \emptyset$. By (\ref{14}), $|\mathcal{H}^{(s)}(T_{\mathcal{B}})| > 0$. Thus, $T_{\mathcal{B}} \in \mathcal{H}^{(t)}$ as $\mathcal{H}$ is hereditary and $T_{\mathcal{B}}$ is a $t$-element subset of every set in $\mathcal{H}^{(s)}(T_{\mathcal{B}})$. Hence  
\begin{equation} |\mathcal{H}^{(s)}(T_{\mathcal{B}})| \leq |\mathcal{H}^{(s)}(T')|. \label{16b}
\end{equation} 
%
%Now, by (\ref{14}) and (\ref{15}),
%
%\[|\mathcal{A}| + |\mathcal{B}| < {s \choose t} |\mathcal{H}^{(r)}(T_D)| + \frac{r(s-t)}{n - s} {r \choose t} |\mathcal{H}^{(s)}(T_{\mathcal{B}})|.\]
%
%By straightforward induction, $2^a \geq 2a$ for every positive integer $a$. We have $n \geq 2^r(r-t)(s-t){r \choose t} + r + s - t \geq 2r(s-t){r \choose t} + s$, so $\frac{r(s-t)}{n - s} {r \choose t} \leq \frac{1}{2}$. Also, $n \geq 2(s-t){s \choose t} + r$ and $t \leq r < s$, so $n - r > s-t > 0$ and ${s \choose t} \leq \frac{1}{2} \left( \frac{n-r}{s-t} \right)$, and hence ${s \choose t} \leq \frac{1}{2} \prod_{i=0}^{s-r-1} \left( \frac{n-r-i}{s-t-i} \right) = \frac{1}{2} \frac{{n-r \choose s-r}}{{s-t \choose s-r}} \leq \frac{1}{2} \frac{{\mu(\mathcal{H})-r \choose s-r}}{{s-t \choose s-r}}$. Thus,
We have
\begin{align}|\mathcal{A}| + |\mathcal{B}| &< {s \choose t} |\mathcal{H}^{(r)}(T_D)| + \frac{r(s-t)}{\mu(\mathcal{H}) - s} {r \choose t} |\mathcal{H}^{(s)}(T_{\mathcal{B}})| \quad \mbox{(by (\ref{14}) and (\ref{15}))} \nonumber \\
&< \frac{1}{2} \frac{{\mu(\mathcal{H}) - r \choose
s-r}}{{s-t \choose s-r}}|\mathcal{H}^{(r)}(T_{D})| + \frac{1}{2} |\mathcal{H}^{(s)}(T_{\mathcal{B}})| \quad \mbox{(by Proposition~\ref{calc} (i) and (ii))} \nonumber \\
&\leq \frac{1}{2}|\mathcal{H}^{(s)}(T')| + \frac{1}{2}|\mathcal{H}^{(s)}(T')| \quad \mbox{(by (\ref{16a}) and (\ref{16b}))} \nonumber \\
&= |\mathcal{H}^{(s)}(T')| < |\mathcal{H}^{(r)}(T')| + |\mathcal{H}^{(s)}(T')|,\nonumber
\end{align}
which is a contradiction since $\emptyset \neq \mathcal{H}^{(r)}(T') \subseteq \mathcal{H}^{(r)}$, $\emptyset \neq \mathcal{H}^{(s)}(T') \subseteq \mathcal{H}^{(s)}$, $\mathcal{H}^{(r)}(T')$ and $\mathcal{H}^{(s)}(T')$ are cross-$t$-intersecting, and $(\mathcal{A}, \mathcal{B}) \in M(\mathcal{H}^{(r)},\mathcal{H}^{(s)},t)$. 

Therefore, $r = s$. Suppose that $\mathcal{B}$ is not a trivial $t$-intersecting family. By Lemma~\ref{transversalemma4}, there exists a $t$-element set $T_{\mathcal{A}}$ such that
\[|\mathcal{A}| < \frac{s(r-t)}{\mu(\mathcal{H})-r}{s \choose t} |\mathcal{H}^{(r)}(T_{\mathcal{A}})|.\]
Thus, $r - t > 0$. Let $T'$ be as defined above (for the case $r < s$). We have
\begin{align}|\mathcal{A}| + |\mathcal{B}| &< \frac{s(r-t)}{\mu(\mathcal{H})-r}{s \choose t} |\mathcal{H}^{(r)}(T_{\mathcal{A}})| + \frac{r(s-t)}{\mu(\mathcal{H}) - s} {r \choose t} |\mathcal{H}^{(s)}(T_{\mathcal{B}})| \nonumber \\
& = \frac{r(r-t)}{\mu(\mathcal{H})-r}{r \choose t} \left( |\mathcal{H}^{(r)}(T_{\mathcal{A}})| + |\mathcal{H}^{(r)}(T_{\mathcal{B}})| \right) \quad \mbox{(as $r=s$)} \nonumber \\
&< \frac{1}{2} \left( |\mathcal{H}^{(r)}(T_{\mathcal{A}})| + |\mathcal{H}^{(r)}(T_{\mathcal{B}})| \right) \quad \mbox{(by Proposition~\ref{calc} (i))} \nonumber \\
&< |\mathcal{H}^{(r)}(T')| + |\mathcal{H}^{(r)}(T')|,\nonumber
\end{align}
which is a contradiction because, as in the case $r < s$ above, $\emptyset \neq \mathcal{H}^{(r)}(T') \subseteq \mathcal{H}^{(r)}$, $\mathcal{H}^{(r)}(T')$ and $\mathcal{H}^{(r)}(T')$ are cross-$t$-intersecting, and $(\mathcal{A}, \mathcal{B}) \in M(\mathcal{H}^{(r)}, \mathcal{H}^{(r)},t)$.

Therefore, $\mathcal{B}$ is a trivial $t$-intersecting family. Thus, since $r = s$, we can apply the argument in Case~1 to obtain that there exists some $I \in \mathcal{H}$ such that $t \leq |I| \leq r$, $\mathcal{B} = \mathcal{H}^{(r)}(I)$, and $\mathcal{A} = \{H \in \mathcal{H}^{(r)} \colon |H \cap I| \geq t\}$. Since $\mathcal{A}$ is not a trivial $t$-intersecting family, $t < |I|$. 
%If $t = |I|$, then $\mathcal{A} = \mathcal{B} = \mathcal{H}^{(r)}(I) = \{H \in \mathcal{H}^{(r)} \colon |H \cap I| \geq t\}$. 
It remains to show that $(\mathcal{A}, \mathcal{B}) \neq (\mathcal{H}^{(r)}(I), \{H \in \mathcal{H}^{(r)} \colon |H \cap I| \geq t\})$ (as the theorem states that the two possibilities resulting from it are mutually exclusive.) 

Since $t < |I|$, $t < r$. Let $T \in {I \choose t}$. Let $B$ be a base of $\mathcal{H}$ such that $I \subseteq B$. Since $\mu(\mathcal{H}) \geq c(r,r,t) \geq r + 2{r \choose t} \geq 3r$, $|B| \geq 3r$. Since $|I| \leq r$, $|B \backslash I| \geq 2r$. Let $X \in {B \backslash I \choose r-t}$. Since $\mathcal{H}$ is hereditary and $T \cup X \subseteq B \in \mathcal{H}$, $T \cup X \in \mathcal{H}$. Thus, $T \cup X \in \mathcal{A} \backslash \mathcal{H}^{(r)}(I)$, and hence $\mathcal{A} \neq \mathcal{H}^{(r)}(I)$. Therefore, $(\mathcal{A}, \mathcal{B}) \neq (\mathcal{H}^{(r)}(I), \{H \in \mathcal{H}^{(r)} \colon |H \cap I| \geq t\})$, as required.~\hfill{$\Box$}

%Then $\mathcal{B} \subseteq \mathcal{H}^{(r)}(I)$ for some set $I$ with $|I| = t$. Let $\mathcal{C} = \mathcal{H}^{(r)}(I)$. Then $\mathcal{A} \cup \mathcal{C}$ and $\mathcal{B}$ are cross-$t$-intersecting. Since $(\mathcal{A}, \mathcal{B}) \in M(\mathcal{H}^{(r)},\mathcal{H}^{(s)},t)$, we have $\mathcal{A} \cup \mathcal{C} = \mathcal{A}$, so $\mathcal{C} \subseteq \mathcal{A}$. Also, $\mathcal{B} \subseteq \mathcal{C}$, so $|\mathcal{B}| \leq |\mathcal{C}| \leq |\mathcal{A}|$. Since $(\mathcal{A}, \mathcal{B}) \in C(\mathcal{H}^{(r)},\mathcal{H}^{(s)},t)$ and $r=s$, we have $|\mathcal{A}| \leq |\mathcal{B}|$, so $|\mathcal{B}| = |\mathcal{C}| = |\mathcal{A}|$. Thus, we have $\mathcal{A} = \mathcal{C} = \mathcal{B}$, which is a contradiction since $\mathcal{A}$ is not a trivial $t$-intersecting family.~\hfill{$\Box$}

%Let $Y := I \cup X^*$ and $Z := S \cup X^*$. Clearly, %$\mathcal{K}^{(p)} = \{H \backslash Y \colon H \in %\mathcal{H}^{(r)}( Y )\}$ and $\mathcal{K}^{(q)} = \{H %\backslash Z \colon H \in \mathcal{H}^{(r)}( Z )\}$. So $|%\mathcal{K}^{(p)}| = |\mathcal{H}^{(r)}( Y )|$ and $|%\mathcal{K}^{(q)}| = |\mathcal{H}^{(r)}( Z )|$. Therefore, %we have

\end{document}